\documentclass{article}
\usepackage[cp1251]{inputenc}
\usepackage[english]{babel}
\usepackage[12pt]{extsizes}
\usepackage{mathrsfs}
\usepackage{amsfonts,amssymb}
\usepackage[tbtags]{amsmath}
\usepackage[T2A]{fontenc}
\usepackage[unicode, pdftex]{hyperref}
\usepackage{mathtools}
\usepackage{enumitem}
\usepackage{graphicx}
\usepackage{caption} 
\usepackage{subcaption}
\usepackage{float,wrapfig}
\usepackage{authblk}

\DeclarePairedDelimiter{\ceil}{\lceil}{\rceil}

\newtheorem{theorem}{Theorem}

\newtheorem{lemma}{Lemma}

\begin{document}

\title{Tight concentration of star saturation number \\ in random graphs}
\author[1]{S. Demyanov}
\author[2]{M. Zhukovskii}
\affil[1]{Moscow Institute of Physics and Technology, Department of Discrete Mathematics}
\affil[2]{The University of Sheffield, Department of Computer Science}

\date{}
\maketitle
\begin{abstract}
   For given graphs $F$ and $G$,  the minimum number of edges in an inclusion-maximal $F$-free subgraph of $G$ is called the $F$-saturation number and denoted $\mathrm{sat}(G, F)$. For the star $F=K_{1,r}$, the asymptotics of $\mathrm{sat}(G(n,p),F)$ is known. We prove a sharper result: whp $\mathrm{sat}(G(n,p), K_{1,r} )$ is concentrated in a set of 2 consecutive points.
   \end{abstract}
\section{Introduction}

The concept of ``saturation number'' was introduced by Zykov \cite{Zyk} and then independently by Erd\H{o}s, Hajnal and Moon \cite{ErHajMoon}. They asked about the minimum number of edges in an $F$-free inclusion-maximal graph on $n$ vertices. In other words, a graph $G$ is {\it $F$-saturated} if $G$ is $F$-free (i.e. does not contain any copy of $F$), but addition of any edge creates a copy of $F$. For example, any complete bipartite graph is a $K_3$-saturated graph. The minimum number of edges in an $F$-saturated graph is called {\it the $F$-saturation number} and is denoted by $\mathrm{sat}(n,F)$.

For example, if $F=K_m$ (i.e. a complete graph on $m$ vertices), then $\mathrm{sat}(n,F)$ is known, this result was obtained by Erd\H{o}s, Hajnal and Moon \cite{ErHajMoon}: for all $n\geq m\geq 2$
\begin{equation*}\label{classic_complete}
	\mathrm{sat}(n,K_m) = (m-2)(n-m+2)+\binom{m-2}{2}=\binom{n}{2}-\binom{n-m+2}{2}.
\end{equation*}

For stars (we denote by $K_{1,r}$ a star with $r$ leaves), the problem was solved by K\'{a}szonyi and Tuza \cite{KasTu}:, 
\begin{equation*}\label{classic_star}
	\mathrm{sat}(n,K_{1,r})=
	\begin{cases}
		\binom{r}{2}+\binom{n-r}{2}, & r+1\leqslant n\leqslant\frac{3r}{2};\\
		\ceil[\big]{\frac{(r-1)n}{2}-\frac{r^2}{8}},& n\geqslant\frac{3r}{2}.
	\end{cases}
\end{equation*}

The notion of saturation can be generalized to arbitrary host graphs. For a given host graph $G$, a spanning subgraph $H$ of $G$ is called {\it $F$-saturated in $G$} if $H$ is $F$-free, but every graph obtained by adding an edge from $E (G) \setminus E(H)$ to $H$, has at least one copy of $F$ as a subgraph. The minimum number of edges in an $F$-saturated subgraph of $G$ is denoted by $\mathrm{sat}(G, F)$. Thus $\mathrm{sat}(n, F)$ = $\mathrm{sat}(K_n, F)$. In \cite {KorSud} Kor\'{a}ndi and Sudakov initiated the study of saturation number of ``typical'' host graphs. We say that a graph property $Q$ holds with high probability (whp), if ${\sf P}\left(G(n,p)\in Q \right) \to 1$ as $n\to\infty$. As usual, we denote by $G(n,p)$ the binomial random graph on $[n]:=\{1,\ldots,n \}$, i.e. a graph with every edge drawn independently with probability $p$.  Kor\'{a}ndi and Sudakov proved that whp $\mathrm{sat} ( G(n,p= \mathrm{const}),K_s)= n \log_{\frac{1}{1-p}} n  (1+o(1))$.

For stars, the saturation number of the random graph was also studied. In \cite{Zito} Zito proved that whp

\begin{equation*}\label{classic_complete}
	    \frac{n}{2}- \log_{\frac{1}{1-p}} \left( np \right) \leqslant \mathrm{sat} \left( G(n,p), K_{1,2} \right) \leqslant 
	\frac{n}{2}- \log_{\frac{1}{1-p}} \left( \sqrt{n} \right) .
\end{equation*}

Note that $\mathrm{sat}(G, K_{1,2})$ is the minimum cardinality of a maximum matching in $G$.

In \cite{MohTay} Mohammadian and Tayfeh--Rezaie  prove that for any fixed $p \in (0,1)$ and any fixed integer $r \geqslant 2$ whp
\begin{equation*}\label{classic_complete}
	\mathrm{sat} \left( G(n,p), K_{1,r} \right)= \frac{(r-1)n}{2}-\left(1+o(1)\right)(r-1)\log_{\frac{1}{1-p}}n.
\end{equation*}

Here we want to emphasize the fact that, for any maximum matching in $G$, the deletion of its vertex set leaves only an independent set of $G$. On the other hand, it is well known that whp any large enough subset of $G(n,p)$ contains a matching (see e.g. \cite[Remark 4.3]{RanGra}) and that the independence number of $G(n,p)$ is concentrated in a set of two consecutive points \cite{Mat1,Mat2, Mat3, GriMcD}: for a fixed $ 0 < p < 1 $ and any $\varepsilon > 0$, whp

\begin{equation}\label{indnum}
\left\lfloor \alpha_p (n) - \varepsilon \right \rfloor \leqslant \alpha (G(n,p)) \leqslant \left\lfloor \alpha_p (n) + \varepsilon \right \rfloor,
\end{equation}
where $\alpha_p (n) := 2 \log_b n - 2\log_b \log_b n + 2\log_b (e/2)+1 $, $b = 1/(1 - p)$ .

Thus, whp $\mathrm{sat} (G(n,p),K_{1,2})$ equals to a half of the size of the complement to a maximum independent set (if this size is odd, then it is 1/2 less), and so it is concentrated in a set of two consecutive points as well. We further consider $r\geqslant 3$.\\

In this paper, we show that for all $r$, whp $\mathrm{sat} (G(n,p),K_{1,r})$ is concentrated in a set of two consecutive points (this is as sharp as possible) and thus significantly improve the result of Mohammadian and Tayfeh--Rezaie \cite{MohTay}. We let 
$$
\varphi_m(k)={n\choose k}{{k\choose 2}\choose m}p^m(1-p)^{{k\choose 2}-m}.
$$
The main result of our paper is stated below.

\begin{theorem}\label{mainth}
Let $p \in (0,1)$, $r \geqslant 3$. Let $\delta>0$, $0<\varepsilon'\ll\varepsilon\ll\delta$ and $n>n_0(\delta)$ be large enough. Let 
$$	
x_0= \left\lfloor \alpha_p (n) + \varepsilon \right \rfloor\quad \text{ and }\quad
r'=\left\lceil\frac{r-3}{2}\right\rceil - I(n-x_0,r-1\text{ are odd}).
$$	
	
If $\varphi_{r'}(x_0)<\varepsilon'$, then, with probability at least $1-\delta$, 
	
		\begin{equation*}\label{maineq_1}
		 \mathrm{sat}(G(n,p), K_{1,r} ) =\left\lceil\frac{(r-1)(n-x_0) }{2}\right\rceil+\mu,\quad \mu:=r'+1.
	\end{equation*}
	
Otherwise, let $\mu\leqslant r'$ be the smallest non-negative integer such that $\varphi_{\mu}(x_0)\geqslant\varepsilon'$. Then, with probability at least $1-\delta$,

		\begin{equation*}\label{maineq_2}
		 \mathrm{sat}(G(n,p), K_{1,r} ) \in\left\{\left\lceil\frac{(r-1)(n-x_0) }{2}\right\rceil+\mu,\left\lceil\frac{(r-1)(n-x_0) }{2}\right\rceil+\mu+1\right\}.%\left\lceil \frac{(r-1)(n-\alpha (G(n,p)) )}{2} \right \rceil  .
	\end{equation*}

\end{theorem}

To prove the theorem, we first show (and this is the trickiest part of the paper) that the almost optimal strategy is\\ (1) to take a set that induces at most $r'$ edges and has the maximum size, and\\ (2) to preserve $r-1$ edges adjacent to each of the vertices outside this set,\\ and, after that, constructively prove the upper bound.\\

The rest of the paper is organised as follows. Section \ref{prelim} contains definitions and theorems needed to prove the main result. In Sections \ref{lbou} and \ref{ubou} the lower and upper bounds are proved respectively.

\section{Preliminaries}\label{prelim}

\subsection{Almost independent sets}
\label{sc:prelim-sets}

Let $\varepsilon>0$ be small enough and let $\xi_m(k)$ be the number of sets of size $k$ that induce exactly $m$ edges in $G(n,p)$. Let $\alpha_m$ be the maximum cardinality of a set of vertices that induces exactly $m$ edges in $G(n,p)$. In particular, $\alpha_0=\alpha(G(n,p))$ is the independence number. Note that $\varphi_m(k)={\sf E}\xi_m(k)$ and that $\alpha_m$ is the maximum $k$ such that $\xi_m(k)\geqslant 1$. 

In~\cite{KSZ} it was proven that, for any constant $m\in\mathbb{Z}_+$, whp $\alpha_m$ is concentrated in a set of 2 consecutive points. Moreover, these points are the same for different values of $m$: for every $m\in\mathbb{Z}_+$, whp
\begin{equation*}
\left\lfloor \alpha_p (n) - \varepsilon \right \rfloor \leqslant \alpha_m \leqslant \left\lfloor \alpha_p (n) + \varepsilon \right \rfloor,
\end{equation*}
More precisely, the following is true.

\begin{theorem}[\cite{KSZ}]
Let $\varepsilon',\varepsilon,\delta,n,x_0$ be defined as in Theorem~\ref{mainth} with the additional requirement that $n\gg m$. We have $\alpha_m\in\{x_0-1,x_0\}$ with probability at least $1-\delta$. If $\varphi_m(x_0)<\varepsilon'$, then $\alpha_m=x_0-1$ with probability at least $1-\delta$. If $\varphi_m(x_0)>1/\varepsilon'$, then $\alpha_m=x_0$ with probability at least $1-\delta$.
\label{th:m-sets}
\end{theorem}

Note that, for every $m$, $\varphi_m(x_0)\ll\varphi_{m+1}(x_0)$. Consider separately two cases distinguished in Theorem~\ref{mainth}.
\begin{enumerate}
\item If $\varphi_{r'}(x_0)<\varepsilon'$, then $\varphi_m(x_0)<\varepsilon'$ as well for all $m<r'$ implying that $\alpha_m=x_0-1$ for all $m\leqslant r'$ with probability at least $1-\delta$. 
\item Otherwise, %let $m\leqslant r'$ be the smallest non-negative integer such that $\varphi_{m}(x_0)\geqslant\varepsilon'$. Then, 
$\alpha_{m}=x_0-1$ for all $m<\mu$, $\alpha_{\mu}\in\{x_0-1,x_0\}$ and $\alpha_{\mu+1}=x_0$ with probability at least $1-\delta$.
\end{enumerate}

% The size of the concentration interval for the independence number of $G(n,p)$ depends on the value of $n$. For some sequences of $n$, the independence number equals $x_0$ whp. For other sequences, whp it belongs to the set $\{x_0-1,x_0\}$, and both values are achievable with probability bounded away from 0. Indeed, in the first case, ${\sf E}\xi_{x_0}\to \infty$ and ${\sf E}\xi_{x_0+1}\to 0$ while $\mathrm{Var}\xi_{x_0}=o(({\sf E}\xi_{k_0})^2)$ implying that ${\sf P}(\xi_{x_0}\geq 1)\to 0$ and ${\sf P}(\xi_{x_0+1}\geq 1)\to 0$ by Chebyshev's and Markov's inequalities respectively (see the details in~\cite[Chapter ?]{RanGra}). Moreover, the relations ${\sf E}\xi_{x_0-1}\to\infty,$ ${\sf E}\xi_{x_0+1}\to 0$ and $\mathrm{Var}\xi_{k}=o(({\sf E}\xi_{k})^2)$ for $k\in\{x_0-1,x_0\}$ hold true in the second case as well, implying that it is only possible when ${\sf E}\xi_{x_0}=\Theta(1)$ since otherwise the independence number equals either $x_0$ whp or $x_0-1$ whp.

%Assume that $(n_i,\,i\in\mathbb{N})$ is a sequence such that ${\sf E}\xi_{x_0}=\Theta(1)$. Let $\xi'_k$ be the number of sets of size $k$ that induce exactly 1 edge in $G(n,p)$. Clearly, ${\sf E}\xi'_k={k\choose 2}\frac{p}{1-p}{\sf E}\xi_k$ implying that ${\sf E}\xi'_{x_0}\to\infty$.

\subsection{Powers of Hamilton cycles}

The $\ell$th power of a graph $H$ is obtained by the addition to $H$ of  edges between all pairs of vertices that are at distance at most $\ell$. A Hamilton $\ell$-cycle in a graph $G$ is the $\ell$th power of a Hamilton cycle in $G$. To prove our main result, we need the following theorem from \cite{KahNar,KO}.
\begin{theorem}[\cite{KahNar,KO}]\label{Ham}
	Let $\ell\geqslant 2$ be fixed. Suppose that  $pn^{1/ \ell}\to\infty$ as $n\to\infty$. Then whp $G(n,p)$ contains a Hamilton $\ell$-cycle.
\end{theorem}
Note that the $\ell$-th power of a Hamilton cycle in a graph on $k(\ell+1)$ admits a $K_{\ell+1}$-factor, i.e. contains a dsjoint union of cliques of size $\ell+1$ covering all vertices of this graph.\\

%According to Theorem \ref{Ham}, for every positive integer $\ell$, whp $G(n,p)$ contains the $\ell$th power of a Hamilton cycle. We need a stronger result:

Let $p=\mathrm{const}\in(0,1)$. Theorem~\ref{Ham} implies the following.
\begin{lemma} \label{lemma2}
Fix a positive integer $\ell$. Whp for every set $W \subset [n]$ of size at most $2\log_b n$ the graph $G(n,p)| _{[n] \backslash W}$ obtained by the deletion of vertices from $W$ contains the $\ell$th power of a Hamilton cycle.
\end{lemma}

The proof is standard and is based on a sequential exposure of edges of $G(n,p)$ sufficiently many times and independently with probability slightly bigger than the threshold probability of an appearance of the $\ell$th power of a Hamilton cycle. For the sake of completeness, we give this argument in Appendix.

\section{ Lower bound} \label{lbou}

Let $G$ be a graph, and $H$ be $K_{1,r}$-saturated in $G$ with the minimum possible number of edges. Obviously, the maximum degree of $H$ does not exceed $r-1$. Let us divide the set of vertices of $H$ into two subsets: $V(H)=V_1\sqcup V_2$, where $V_1$ comprises all the vertices of $H$ with degrees at most $r-2$, and all the vertices in $V_2$ have degree exactly $r-1$. Let the induced subgraph $H_1:=H[V_1]$ have $k$ vertices and $m$ edges. Then the graph $H$ has at least $\frac{r-1}{2}(n-k)+m$ edges, i.e.
\begin{equation}\label{est}
		\mathrm{sat}(G, K_{1,r} ) \geqslant \left\lceil\frac{r-1}{2}(n-k)\right\rceil+m .
	\end{equation}

It is clear that $H_1$ is also the induced subgraph of $G$ itself. If this were not the case, then we could draw an edge inside $H_1$ that does not create $K_{1,r}$  --- a contradiction. So, $|E(G[V_1])|=m$. We will show that if $G=G(n,p)$, then whp it is optimal to take $V_1$ such that $m\leq r'+1$. For that, let us rewrite \eqref{est} as follows:
 
\begin{equation}\label{est2}
		\mathrm{sat}(G, K_{1,r} ) \geqslant \frac{r-1}{2}(n-k)+m = \frac{r-1}{2}(n-x_0)+\left(m-\frac{r-1}{2}(k-x_0) \right).
	\end{equation}

 The term $\frac{r-1}{2}(n-x_0)$ refers to the number of edges in $H$, when $G=G(n,p)$, $H_1$ is an independent set of a maximum size, and there are no edges between $V_1$ and $V_2$ in $H$. We show that the second summand in the right hand side of \eqref{est2} is whp at least $\mu$, where $\mu$ is defined in the statement of Theorem~\ref{mainth}. In other words, whp it is impossible to enlarge a maximum set that induces exactly $\mu$ edges and get an induced subgraph with the number edges less then the increase in its number vertices times $\frac{r-1}{2}$. 

\begin{lemma}\label{lemma1}
Let $p \in (0,1)$ be a constant. Then whp for any $k\geqslant x_0+1$ there are no induced subgraphs on $k$ vertices in $G(n,p)$ with the number of edges less than $\frac{r-1}{2}(k-x_0)+\mu$.
\end{lemma}

Assume that Lemma \ref{lemma1} is true and consider separately three cases.

\begin{enumerate}

\item If $k\leqslant x_0-1$, then (\ref{est2}) implies $\mathrm{sat}(G, K_{1,r} ) \geqslant \left\lceil\frac{r-1}{2}(n-x_0)+\frac{r-1}{2}
\right\rceil\geqslant \left\lceil\frac{r-1}{2}(n-x_0)\right\rceil+\mu$ as needed.

\item If $k=x_0$, then, due to the definition of $\mu$, whp $m\geq\mu$ --- see the discussion in Section~\ref{sc:prelim-sets} after the statement of Theorem~\ref{th:m-sets}. It readily implies the desired inequality $\mathrm{sat}(G, K_{1,r} ) \geqslant \frac{r-1}{2}(n-x_0)+\mu$.

\item Finally, if $k\geqslant x_0+1$, then, due to Lemma~\ref{lemma1}, whp $m-\frac{r-1}{2}(k-x_0)\geq\mu$ and therefore $\mathrm{sat}(G, K_{1,r} ) \geqslant \frac{r-1}{2}(n-x_0)+\mu$ as well.

%\begin{equation*}\label{est3}
%		\mathrm{sat}(G(n,p), K_{1,r} ) \geqslant  \frac{r-1}{2}(n-x_0).
%	\end{equation*}
\end{enumerate}

This finishes the proof of the lower bound in Theorem \ref{mainth}.  In Section \ref{PoL1} we give the proof of Lemma~\ref{lemma1}. % First of all note that $\mu\leq r'+1\leq\left\lceil\frac{r-1}{2}\right\rceil$.

\subsection{ Proof of Lemma \ref{lemma1}} \label{PoL1}
Denote $x_1:=x_0+1$. Let $X_k$ be a random variable equal to the number of induced subgraphs on $k$ vertices with the number of edges fewer than than $\frac{r-1}{2}(k-x_0)+\mu$, $X = \sum\limits_{k \geqslant x_1} X_k $. Notice that Lemma \ref{lemma1} states that $X=0$ whp. We will prove it be Markov's inequality. Thus, it is sufficient to show that $\mathbb{E}X \to 0$ as $n \to \infty$.

Due to the linearity of the expectation
\begin{equation}\label{exp}
    \mathbb {E}X = \sum\limits_{k \geqslant x_1} \mathbb {E}X_k = \sum\limits_{k \geqslant x_1} \binom{n}{k} \sum\limits_{m <\mu+ \frac{r-1}{2}(k-x_0)} \binom{\binom{k}{2}}{m} p^m (1-p)^{\binom{k}{2}-m}.
\end{equation}

Let us show that the function $f(m):= \binom{\binom{k}{2}}{m} p^m (1-p)^{\binom{k}{2}-m}$ increases in our range. To do this, consider the ratio $\frac{f(m+1)}{f(m)}$ and show that it is greater than $1$:
\begin{gather*}
    \frac{f(m+1)}{f(m)} =  \frac{p}{1-p}\times \frac{\binom{k}{2}-m}{m+1} \geqslant \frac{p}{1-p}\times\frac{{k\choose 2}-\frac{r-1}{2}(k-x_0)-\mu}{\frac{r-1}{2}(k-x_0)+\mu+1}>1
   \end{gather*}
for $n$ large enough.

 Let us rewrite \eqref{exp} as follows:
\begin{equation}\label{exp2}
    \mathbb {E}X =  \binom{n}{x_1} \sum\limits_{m < \mu+\frac{r-1}{2}} \binom{\binom{x_1}{2}}{m} p^m (1-p)^{\binom{x_1}{2}-m}+ \sum\limits_{k \geqslant x_1 + 1} \binom{n}{k} \sum\limits_{m < \mu + \frac{r-1}{2}(k-x_0)} \binom{\binom{k}{2}}{m} p^m (1-p)^{\binom{k}{2}-m}.
\end{equation}

Consider separately the first and second terms in \eqref{exp2} and show that they tend to zero. We start with the first summand. Let us recall that
\begin{gather*}
   x_1 =\left\lfloor 2 \log_b n - 2\log_b \log_b n + 2\log_b (e/2)+1+ \varepsilon \right \rfloor +1.
\end{gather*}

Therefore,
\begin{gather*}
    \frac{x_1-1}{2} \geqslant \log_b n - \log_b \log_b n + \log_b (e/2)+ \frac{\varepsilon}{2} .
\end{gather*}

We get
\begin{align*}
       \binom{n}{x_1}(1-p)^{\binom{x_1}{2}} & \leqslant \left( \frac{en}{x_1} \right)^{x_1}(1-p)^{\binom{x_1}{2}} =  \exp \left[ x_1 \left( \ln n + 1 - \ln x_1 - \frac{x_1 - 1}{2} \ln b \right)  \right]  \\ 
       & \leqslant    \exp \left[ x_1 \left( 1 - \ln 2 - \log_b (e/2)\ln b - \frac{\varepsilon}{2} +o(1) \right)  \right] \\
    &= \exp \left[  - x_1\left(\frac{\varepsilon}{2}+ o(1) \right)    \right]= \exp \left[  - \varepsilon \log_b n (1 + o (1))    \right].
\end{align*}

Moreover,
\begin{gather*}
    \sum\limits_{m < \mu + \frac{r-1}{2}} 
    \binom{\binom{x_1}{2}}{m} \left( \frac{p}{1-p} \right)^m 
    \leqslant \left[\mu+\frac{r-1}{2}\right]\max_m \left( \frac{ \binom{x_1}{2} e}{m} \right)^m \left(\frac{p}{1-p}\right)^m  =o(x_1^{2\mu+r}).
\end{gather*}

From the above, the first summand in \eqref{exp2} is bounded from above as follows:
\begin{equation}\label{exp5}
     \binom{n}{x_1}(1-p)^{\binom{x_1}{2}} \sum\limits_{m < \mu + \frac{r-1}{2}} \binom{\binom{x_1}{2}}{m} \left(\frac{p}{1-p}\right)^m  \leqslant \exp[-\varepsilon\log_b n(1+o(1))+(r+2\mu)\ln x_1]=o(1)
\end{equation}
as needed.

Now let us switch to the second summand in \eqref{exp2}. Let $k \geqslant x_1 + 1$. Then

\begin{align*} 
\binom{n}{k} & \sum\limits_{m < \mu + \frac{r-1}{2}(k-x_0)} \binom{\binom{k}{2}}{m} p^m (1-p)^{\binom{k}{2}-m}  \\
&\leqslant \left(\frac{ne}{k} \right)^k\left[\mu+ \frac{r-1}{2}(k-x_0)\right] \left(\frac{\binom{k}{2}e}
{\mu+\left\lfloor \frac{r-1}{2}(k-x_0) \right\rfloor}\right)^{\mu+\left\lfloor \frac{r-1}{2}(k-x_0) \right\rfloor} \left( \frac{p}{1-p}\right)^{\mu+\left\lfloor \frac{r-1}{2}(k-x_0) \right\rfloor}(1-p)^{\binom{k}{2}}\\
&=\left[\left( \frac{ne}{x_1}\right)^{x_1}(1-p)^{\binom{x_1}{2}}\right](ne)^{k-x_1} \frac{x_1^{x_1}}{k^k}\left[\mu+ \frac{r-1}{2}(k-x_0)\right]
\times \\
&\hspace{4cm}\times \left(\frac{\binom{k}{2}e}
{\mu+\left\lfloor \frac{r-1}{2}(k-x_0) \right\rfloor}\right)^{\mu+\left\lfloor \frac{r-1}{2}(k-x_0) \right\rfloor}
\left( \frac{p}{1-p}\right)^{\mu+\left\lfloor \frac{r-1}{2}(k-x_0) \right\rfloor}
(1-p)^{\binom{k}{2}-\binom{x_1}{2}}.
\end{align*}

Note that the factor $ \left( \frac{ne}{x_1}\right)^{x_1}(1-p)^{\binom{x_1}{2}}$ does not exceed 1, because this is an upper bound for ${\sf E}\xi_0(x_1)$ (the expected number of independent sets of size $x_1$), and $x_1$ is chosen exactly in a way such that this bound approaches 0 (see \cite[Remark 7.3]{RanGra}). Then the last expression can be estimated from above for large enough $n$ as

\begin{multline} \label{last}
 \frac{ (ne)^{k-x_1} x_1^{x_1}}{k^k}\left[\mu+ \frac{r-1}{2}(k-x_0)\right]\left(\frac{\binom{k}{2}e}
{\mu+\left\lfloor \frac{r-1}{2}(k-x_0) \right\rfloor}\right)^{\mu+\left\lfloor \frac{r-1}{2}(k-x_0) \right\rfloor} \times\\
\times \left( \frac{p}{1-p}\right)^{\mu+\left\lfloor \frac{r-1}{2}(k-x_0) \right\rfloor}(1-p)^{\binom{k}{2}-\binom{x_1}{2}} =: f(k).
\end{multline}

Then

\begin{align*}     
\ln f(k) &  \sim     (k -x_1)\ln n  - (k-x_1)\ln k + \left[\mu+\frac{r-1}{2}(k-x_0)\right] \ln \binom{k}{2} -\\    & \quad \; \; \; \; \,  -  \left[\mu+\frac{r-1}{2}(k-x_0)\right] \ln \left(\mu+\frac{r-1}{2} (k-x_0)\right)- \left( \binom{k}{2} -  \binom{x_1}{2}  \right) \ln \frac{1}{1-p} \\   &\leqslant (k-x_1)\ln n  + \left[\mu+\frac{r-1}{2}(k-x_0)\right] \ln \binom{k}{2}- \left( \binom{k}{2} -  \binom{x_1}{2}  \right) \ln \frac{1}{1-p}  \\   &\sim (k-x_1)\ln n  + (r-1)(k-x_1) \ln k-\frac{(k-x_1)(k+x_1)}{2} \ln \frac{1}{1-p}. \end{align*}

 Denote $c:=k-x_1$, then 
\begin{equation*} 
   \ln f(k) \leqslant c\left(\ln n+ (r-1) \ln k - \frac{2x_1 + c}{2} \ln \frac{1}{1-p}\right)(1+o(1)).
\end{equation*}

Differentiating the function inside the brackets with respect to $k$, we get $ \left( \frac{r-1}{k}- \frac{1}{2}\ln \frac{1}{1-p} \right )$. Since $k \to \infty $, then $  \frac{r-1}{k}- \frac{1}{2}\ln \frac{1}{1-p}  < 0$ for sufficiently large $n$, and hence the function itself is decreasing when $k \geq x_1$. And since it decreases, then its maximum value is reached at the smallest possible $k=x_1$. Note that $\ln x_1=O(\ln\ln n)$, and it does note affect asymptotics. Summing up, we get that
 \begin{equation*} 
   \ln f(k)\leq (k-x_1) \left(\ln n-x_1\ln\frac{1}{1-p}\right)(1+o(1) )\sim -\ln n(k-x_1)(1+o(1)).
\end{equation*}

Therefore, the second summand in \eqref{exp2} is bounded from above as follows:
\begin{equation}\label{es3}
    \sum\limits_{k \geqslant x_1 +1} \binom{n}{k} \sum\limits_{m < \mu+ \frac{r-1}{2}(k-x_0)} \binom{\binom{k}{2}}{m} p^m (1-p)^{\binom{k}{2}-m} \leqslant \sum\limits_{c \geqslant 1} e^{-c\ln n(1+o(1))} \to 0.
\end{equation}

Due to \eqref{exp5} and \eqref{es3},  $\mathbb {E}X \to 0$ as $n\to\infty$ as needed.

\section{Upper bound}\label{ubou}

Assume first that $\mu\leq m$. %If in $G:=G(n,p)$ there exists a set with $\mu$ edges and $x_0$ vertices, then we let $V_1$ to be such a set (note that whp $x_0$ is the maximum possible size of a set with $\mu$ edges due to Theorem~\ref{th:m-sets}). 
 Then whp there exists a set with $\mu+1$ edges and $x_0$ vertices (see the discussion after the statement Theorem~\ref{th:m-sets}), and we let $V_1$ to be such a set. If $\mu=m+1$, then $V_1$ is an independent set on $x_0-1$ vertices (it exists whp due to (\ref{indnum})). Let us show that whp there exists a subgraph in $G$ such that each vertex from $V_2:=[n]\setminus V_1$ has degree exactly $r-1$, and there is at most 1 edge between $V_1$ and $V_2$. Clearly, such a subgraph is saturated in $G$.

According to Lemma \ref{lemma2}, whp $G|_{V_2}$ contains the $(2r-2)$th power of a Hamilton cycle. Let us preserve in this cycle a disjoint union of cliques $K_r$ of size $r$ and a clique $K^*$ with at least $r$ and at most $2r-1$ vertices such that this union covers all vertices of $V_2$.

It remains to show that we may turn $K^*$ into an $(r-1)$-regular graph by deleting some edges of $K^*$ and drawing at most 1 edge from $K^*$ to $V_1$.

 Let us first remove some edges from $K^{*}$ so that only a simple cycle containing all its vertices remains.

If $r-1=2s$ is even, then join each vertex in this cycle with $s$ nearest neighbors. We get $(r-1)$-regular graph.
If $r-1=2s+1$ is odd and $|V(K^{*})|$ is even then join each vertex in the cycle with $s$ nearest neighbors and also with the opposite vertex.
And finally, if $r-1=2s+1$ is odd and $|V(K^{*})|$ is odd as well, then we arbitrarily choose a single vertex $v_1$ in the cycle. Note that $V_1$ is either an independent set with maximum possible size or a set with $\mu$ edges with maximum possible size, and therefore $v_1$ has a neighbor in $V_1$. Draw a single edge of $G$ from $v_1$ to $V_1$. Inside the cycle, we join each vertex with its $s$ nearest neighbors. Then, eventually, consider the cyclic order on $V(K^*)\setminus \{v_1\}$ which is exactly the order induced by the initial cycle from wich we exclude the vertex $v_1$. Draw from every vertex in $V(K^*)\setminus\{v_1\}$ the edge to the opposite vertex in this order.

The desired saturated subgraph of $G$ is constructed, it has exactly $\left\lceil \frac{(r-1)(n-x_0)}{2}  \right \rceil+\mu+1$ edges if $\mu\leq m$ and $\left\lceil \frac{(r-1)(n-x_0)}{2}  \right \rceil+\mu$ edges if $\mu=m+1$ implying the upper bound in Theorem~\ref{mainth}. %that whp

%\begin{equation*}
%    \mathrm{sat}(G(n,p), K_{1,r} ) \leqslant  \left\lceil \frac{(r-1)(n-x_0)}{2}  \right \rceil+\mu.
%\end{equation*}

%\section{ Sharp concentration} \label{concl}

%Note that the size of the concentration interval for independence number depends on the value of $n$. For some sequences of $n$, the independence number equals to a specific $k=f(n)$ whp. For other sequences, whp it belongs to a $\{k,k+1\}$, and both values are achievable with probability bounded away from 0. In the first case, whp $\alpha(G(n,p))=x_0$, and so the upper bound and the lower bound in Theorem \ref{mainth} differ by no more than 1/2 implying that the saturation number equals to the upper bound whp. 
%In the second case, it is easy to see that in the statement of Lemma \ref{lemma1} $x_0$ can be replaced with $x_1$. Indeed, the only difference in the proof is the place where we omit the factor $ \left( \frac{ne}{x_1}\right)^{x_1}(1-p)^{\binom{x_1}{2}}$ in Equation \ref{last}. Although, this factor does not approach 0 anymore, it is still sufficiently small, and does not affect significantly the value of $f(k)$. It implies that whp the saturation number belongs to a set of two values, but they are not consecutive --- the difference between them is $\left\lceil(r-1)/2\right \rceil$. So, the length of the concentration interval given in Theorem \ref{mainth} is tight.

\section*{Appendix: proof of Lemma \ref{lemma2}} \label{app}    
Let $Q$ be the property that a graph contains the $\ell$th power of a Hamilton cycle. According to Theorem \ref{Ham}, the threshold probability of this property is $\hat{p}=n^{-\frac{1}{\ell}}$. Let $p_0 = n^{-\frac{1}{\ell}} \ln{n}$. Then $ {\sf P} \left( G(n,p_0) \notin Q \right) \leqslant \frac{1}{2}$ for large enough $n$.

Let $x$ be the maximum integer such that $1-p \leq (1-p_0)^x$. Let us take logarithm of the both sides of this inequality and find the asymptotic behavior of $x$:
\begin{equation} \label{xin}
    (x+O(1)) \ln(1-p_0) = \ln(1-p).
\end{equation}

Since $p_0 \to 0$ for $n \to \infty$, then $\ln(1-p_0) \sim -p_0$ and \eqref{xin} can be written as
\begin{equation*} 
    x  \sim \frac{1}{p_0} \ln\left(\frac{1}{1-p}\right)= \frac{n^{\frac{1}{\ell}}}{\ln n} \ln \left(\frac{1}{1-p}\right).
\end{equation*}

Consider a union of $x$ independent copies $G_1,\ldots,G_x$ of $G(n,p_0)$. Clearly, there exists a coupling such that this union $G$ is a subgraph of $G(n,p)$. Therefore, if $G(n,p)$ does not contain the $\ell$th power of a Hamilton cycle, then $G$ does not contain it as well, and the same applies to each of $G_1,\ldots,G_k$. Then
\begin{equation*}
    {\sf P} ( G(n,p) \notin Q ) \leqslant 
    \left( {\sf P} ( G(n,p_0) \notin Q ) \right)^x  \leqslant \left( \frac{1}{2} \right)^{n^{1 / \ell}  (\ln b+o(1)) / \ln n}.
\end{equation*}

By the union bound, the probability that there exists a set $W$ of size at most $2\log_b n$ such that $G(n,p)|_{[n]\setminus W}$ does not contain the $\ell$th power of a Hamilton cycle is at most
 \begin{equation*}
        2\log_{b} n \;   \; n^{2\log_{b} n} \; {\sf P}    \left( G(n(1-o (1)),p) \notin Q \right) \leqslant 
        2\log_{b} n \;   \; n^{2\log_{b} n} \; e^{- \ln 2 \ln b n^{1 / \ell}(\ln n)^{-1}(1+o(1))} \to 0, \quad n\to\infty
\end{equation*}
as needed.


\begin{thebibliography}{99}
    
	
\bibitem{ErHajMoon} P. Erd\H{o}s, A. Hajnal, J. Moon, A problem in graph theory, {\it Amer. Math. Monthly} 71 (1964) 1107--1110.

\bibitem{GriMcD}	G. R. Grimmett, C. J. H. McDiarmid, On colouring random graphs, {\it Math. Proc. Cambridge Philos. Soc.} 77 (1975) 313--324.
        
\bibitem{RanGra} S. Janson, T. \L uczak, A. Ruci\'{n}ski, Random Graphs, Wiley, New York (2000).
	
\bibitem{KahNar} J. Kahn, B. Narayanan, J. Park, The threshold for the square of a Hamilton cycle, {\it Proc. Amer. Math. Soc.} 149 (2021) 3201--3208.

\bibitem{KSZ} D. Kamaldinov, A. Skorkin, M. Zhukovskii, Maximum sparse induced subgraphs of the binomial random graph with given number of edges, {\it Discrete mathematics} 344 (2021) 112205.
	
\bibitem{KasTu} L. K\'{a}szonyi, Zs. Tuza, Saturated graphs with minimal number of edges, {\it J. Graph Theory} 10 (1986) 203--210.

\bibitem{KorSud} D. Kor\'{a}ndi, B. Sudakov, Saturation in random graphs, {\it Random Structures \& Algorithms} 51:1 (2016) 169--181. 

\bibitem{KO} D. K\"{u}hn, D. Osthus, On P\'{o}sas conjecture for random graphs, {\it SIAM Journal on	Discrete Mathematics}, 26 (2012) 1440--1457. 

\bibitem{Mat1} D. W. Matula, On the complete subgraphs of a random graph, In {\it Proceedings of the 2nd Chapel Hill Conference on Combinatorial Mathematics and its Applications} (Chapel Hill, N. C., 1970) 356--369.
	
\bibitem{Mat2} D. W. Matula, The employee party problem, {\it Notices AMS} 19:2 (1972) A--382.
	
\bibitem{Mat3} D. W. Matula, The largest clique size in a random graph. Tech. Rep. 1987, Department of Computer Science, Southern Methodist University, Dallas, Texas (1976).

\bibitem{MohTay} A. Mohammadian, B. Tayfeh--Rezaie, Star saturation number of random graphs, {\it Discrete Mathematics} 341:4 (2018) 1166--1170. 

\bibitem{Zito} M. Zito, Small maximal matchings in random graphs, {\it Theoret. Comput. Sci.} 297 (2003) 487--507.

\bibitem{Zyk} A. A. Zykov, On some properties of linear complexes, {\it Mat. Sbornik} (Russian) 24:66 (1949) 163--188.

    
\end{thebibliography}
\end{document}